\documentclass[11pt,reqno]{amsart}
\usepackage{mathrsfs}
\usepackage{amsmath}
\usepackage{amsxtra}
\usepackage{amscd}
\usepackage{amsthm}
\usepackage{amsfonts}
\usepackage{amssymb}
\usepackage{pstricks}
\usepackage{pstricks,pst-node}

\newcommand{\mysection}[1]{\section{#1}\setcounter{equation}{0}}

\newtheorem{theorem}{Theorem}[section]

\theoremstyle{definition}

\theoremstyle{remark}

\numberwithin{equation}{section}

\def\bS{{\mathbb S}}

\def\P{\mathcal P}

\newcommand{\kg}{/\kern-0.55em g}
\newcommand{\kS}{{/\kern-0.55em \mathbb{S}} _\Sigma}
\newcommand{\kC}{/\kern-0.65em \nabla}
\newcommand{\kD}{/\kern-0.65em D}
\newcommand{\kl}{/\kern-0.50em \lambda_1}
\newcommand{\kgamma}{/\kern-0.65em \gamma}
\newcommand{\n}{\nabla}
\newcommand{\ns}{\nabla^{\Sigma}}
\newcommand{\Ds}{D^{\Sigma}}
\newcommand{\gs}{g\!\!\!/\,}
\newcommand{\abs}[1]{\left|#1\right|}
\newcommand{\obS}{\overline{\bS}}
\newcommand{\bg}{\bar{g}}
\newcommand{\oconn}{\overline{\nabla}}
\newcommand{\psb}[2]{\left\langle{#1},{#2}\right\rangle_{\bar{g}}}

\newcommand{\ps}[2]{\left\langle{#1},{#2}\right\rangle}

\newcommand{\ls}{\setlength{\baselineskip}{12pt}
                 \setlength{\parskip}{3mm}}

\begin{document}


\title[Dirac operator]{Eigenvalue Estimate of the Dirac operator and Rigidity of Poincare-Einstein Metrics}

\author[]{Daguang Chen$^{\dag}$, Fang Wang$^{\ddag}$ and Xiao Zhang$^{\flat}$}

\address[]{$^{\dag}$ Department of Mathematical Sciences, Tsinghua University, Beijing, 100084, P.R. China }

\email{dgchen@math.tsinghua.edu.cn}

\address[]{$^{\ddag}$ School of Mathematical Sciences, Shanghai Jiao Tong University, Shanghai 200240, P.R. China}

\email{fangwang1984@sjtu.edu.cn}

\address[]{$^{\flat}$ Institute of Mathematics,
Academy of Mathematics and Systems Science, Chinese Academy of Sciences, Beijing 100190, and
School of Mathematical Sciences, University of Chinese Academy of Sciences,
Beijing 100049, P.R. China}

\email{xzhang@amss.ac.cn}

\begin{abstract}

We re-visit the eigenvalue estimate of the Dirac operator on spin manifolds with boundary in terms of the first eigenvalues
of conformal Laplace
 operator as well as the conformal mean curvature operator. These problems were studied earlier by
Hijazi-Montiel-Zhang and Raulot and we re-prove them under weaker assumption that a boundary chirality operator exists.
Moreover, on these spin manifolds with boundary, we show that any $C^{3,\alpha}$ conformal compactification of some
Poincare-Einstein metric must be the standard hemisphere when the first nonzero eigenvalue of the Dirac operator achieves its lowest
value, and any $C^{3,\alpha}$ conformal compactification of some Poincare-Einstein metric must be the flat ball in Euclidean space
when the first positive eigenvalue of the boundary Dirac operator achieves certain value relating to the second Yamabe invariant.
In two cases the Poincare-Einstein metrics are standard hyperbolic metric.

\end{abstract}

\keywords{Eigenvalue, Dirac operator, boundary condition, Yamabe invariant, Poincare-Einstein metric}

\maketitle

%
%
\mysection{Introduction}\ls
Let $(M,g)$ be a closed (compact without boundary) $n$-dimensional Riemannian spin manifold with the positive scalar curvature $R>0$. Let $\lambda(D)$ be the eigenvalue of the Dirac operator $D$.  In 1963, Lichnerowicz \cite{L} firstly proved
\begin{equation}\label{Lich}
\lambda ^2(D)>\frac 14 \underset{M}\inf\, R .
\end{equation}
By modifying the Riemammian spin connection suitably, Friedrich \cite{Fr80} improved the Lichnerowicz inequality \eqref{Lich} and obtained the sharp estimate
\begin{equation}\label{Fri}
\lambda  ^2(D)\geq \frac{n}{4(n-1)}\underset{M} \inf R .
\end{equation}
If the equality holds in \eqref{Fri}, the manifold is Einstein. In 1986, using conformal covariance of the Dirac operator, Hijazi \cite{H86} showed, for $n\geq3$,
\begin{equation}\label{Hijazi}
\lambda ^2(D)\geq\frac {n}{4(n-1)} \mu_1,
\end{equation}
where $\mu_1$ is the first eigenvalue of the conformal Laplace operator. If the equality holds in \eqref{Hijazi}, there exists the real Killing spinor and the manifold becomes Einstein.

For any $n$-dimensional ($n \geq 3$) compact manifold $M$ with boundary $\Sigma$, let $\{e_\kappa\}_{\kappa =1}^n$ be the local orthonormal frame along $\Sigma$ such that $e_n$ is a global outward normal to $\Sigma$ and $\{e_i\}_{i=1}^{n-1}$ is tangent to $\Sigma$. We denote by $\kC$ the Levi-Civit\`{a} connection with respect to the induced metric $\kg$ on the hypersurface $\Sigma$. The Gauss formula gives
\begin{equation}\label{Gauss-M}
\n_i e_j=\kC _i e_j-\,h_{ij}e_n
\end{equation}
where $h_{ij}$ is the second fundamental form of $\Sigma$ defined by
\begin{equation*}
h_{ij}=g(\n_ie_n,e_j)=-g(\n_ie_j,e_n).
\end{equation*}
The mean curvature $H$ of hypersurface $\Sigma$ is given by
\begin{equation*}
H=\frac{1}{n-1}tr _g h .	
\end{equation*}
The conformal Laplace operator $L$ and the conformal mean curvature operator $B$ are defined as
\begin{eqnarray*}
&&  L=-\frac{4(n-1)}{n-2}\Delta +R, \label{conf-Laplace}\\
&&  B=\frac{2}{n-2}e_n+H. \label{conf-mOperator}
\end{eqnarray*}
The variational characterizations of the first eigenvalue of $L$ and $B$ are given by
\begin{eqnarray*}
&& \mu _1 (L) =\underset{f\in C^1(\overline{M}),f\neq 0}{\inf} \frac{\int_{M}\Big( \frac{4(n-1)}{n-2}|\nabla f|^2
	+R f^2\Big)d\mu_g +2(n-1)\int_{\Sigma} Hf^2 d\sigma_{\gs}}{\int_{M} f^2 d\mu_g}, \label{EigenV-CL}\\
&& \nu _1(B) =\underset{f\in C^1(\overline{M}),f\neq 0}{\inf} \frac{\int_{M}\Big( \frac{2}{n-2}|\nabla f|^2
	+\frac{1}{2(n-1)}R f^2\Big)d\mu_g +\int_{\Sigma} Hf^2 d\sigma_{\gs} }{\int_{\Sigma}f^2d\sigma_{\gs}}, \label{EigenV-CM}
\end{eqnarray*}
respectively. In \cite{Escobar92JDG, Escobar92}, Escobar proved the first positive eigenfunctions always exist, i.e. there exists a unique $f>0$ satisfying
\begin{equation}\label{Eigen-bdry0}
\left \{ \aligned
L f = & \,\mu_1 (L)f \qquad          \text{on} \quad  M, \\
B f = & \,0       \qquad\qquad \, \text{on} \quad \Sigma,
\endaligned
\right.
\end{equation}
and there exists a unique $f>0$ satisfying
\begin{equation}\label{Eigen-bdry1}
\left \{ \aligned
L f = & \,0\qquad\qquad \, \text{on} \quad M, \\
B f = & \, \nu_1 (B)f \qquad    \text{on}     \quad \Sigma,
\endaligned
\right.
\end{equation}
provided $\nu_1(B)>-\infty$. (It was first pointed out by Jin that $\nu_1(B)$ could be $-\infty$, and this is the case to remove a small geodesic ball on a compact manifolds without boundary with negative scalar curvature \cite{Escobar92}.)

For compact manifold $M$ with boundary $\Sigma$, the (normalized) first and the second Yamabe invariants are given by
 \begin{eqnarray*}
&& Y(M,\Sigma)=\underset{f\in C^1(\overline{M}),f\neq 0}{\inf} \frac{\int_{M}\left( \frac{4(n-1)}{n-2}|\nabla f|^2
		+Rf^2\right)d\mu_g +2(n-1)\int_{\Sigma} Hf^2 d\sigma_{\gs }} {\left(\int_{M} f^\frac{2n}{n-2}d\mu_g \right)^\frac{n-2}{n} },\label{Yamabe1}\\
&& Q(M,\Sigma)=\underset{f\in C^1(\overline{M}),f\neq 0}{\inf} \frac{\int_{M}\left( \frac{2}{n-2}|\nabla f|^2
	+\frac{1}{2(n-1)}R f^2\right)d\mu_g +\int_{\Sigma} Hf^2 d\sigma_{\gs } }{\left(\int_{\Sigma}f^\frac{2(n-1)}{n-2}d\sigma_{\gs }\right)^\frac{n-2}{n-1} }, \label{Yamabe2}
\end{eqnarray*}
respectively. If $\mu_1(L)\geq 0$, $\nu _1(B)\geq 0$, by H\"older inequality, we have
\begin{eqnarray}
&& \mu_1(L) \geq \frac{Y(M, \Sigma)}{\mbox{Vol}(M) ^\frac{2}{n} },\label{Yamabe-Eigen1}\\
&& \nu_1(B) \geq \frac{Q(M, \Sigma)}{\mbox{Vol}(\Sigma) ^\frac{1}{n-1} }.\label{Yamabe-Eigen2}
\end{eqnarray}	
Equality occurs in \eqref{Yamabe-Eigen1} if and only if the corresponding eigenfunction is constant in $M$ and equality occurs in
\eqref{Yamabe-Eigen2} if and only if the corresponding eigenfunction is constant on $\Sigma$.	

For compact spin manifold $M$ with boundary $\Sigma$,  suitable boundary conditions should be imposed in order to make the Dirac operator self-adjoint and elliptic. There exist two basic types of boundary conditions, the global Atiyah-Patodi-Singer (APS) boundary condition and the local boundary condition \cite{APS1,APS2,APS3,GHHP,FS}. The Friedrich inequality was generalized to spin manifolds with boundary under the two types of boundary conditions as well as certain mixed boundary condition \cite{HMZ01CMP,HMZ02}. For conformal aspect of the Dirac operator on manifolds with boundary, the APS boundary condition is not conformal invariant, but the local boundary condition can be used to generalize the Hijazi inequality to spin manifolds with boundary for $n\geq 3$ with a boundary chirality operator \cite{HMZ01CMP} as well as for $n\geq 2$ with a chirality operator \cite{R06}
\begin{equation}\label{HMZ01}
\lambda _1 ^2(D) \geq
\left \{ \aligned
\frac{n}{4(n-1)}\mu _1 (L),\qquad n \geq 3, \\
\frac{2\pi}{Area(M^2, g)}, \qquad n=2,
\endaligned
\right.
\end{equation}
when $\mu_1(L)>0$. Moreover, for internal boundary $\Sigma $ of compact domain in a spin manifold, the conformal integral Schr\"{o}dinger-Lichnerowicz formula and local boundary condition also yield
\begin{equation}\label{HMZ02}
\kl(\Ds) \geq \frac{n-1}{2}\nu _1 (B)
\end{equation}
when $\nu _1 (B)>0$, where $\kl(\Ds)$ is the first positive eigenvalue of the Dirac operator of $\Sigma$ \cite{HMZ02}. It was assumed that $\Sigma$ is an internal hypersurface in order to use the unique continuation property of the Dirac operator. But this property does not seem to be verified when $\Sigma$ is the boundary of $M$ and the Riemannian structure and spin structure are not products near $\Sigma$ (c.f. Remark 8.4 in \cite{BW}).

In this paper, we re-visit and prove \eqref{HMZ01} and \eqref{HMZ02} when $\Sigma$ equips with a boundary chirality operator. For $n\geq 3$, we also study the rigidity of $(M,\Sigma,g)$ as a $C^{3, \alpha}$ conformal compactification of Poincar\'{e}-Einstein manifold $(\mathring{M},  g_+)$:
$$
Ric_{g_+}=-(n-1)g_+,\quad\textrm{in $\mathring{M}$}
$$
and $g=\rho^2 g_+$ can be $C^{3, \alpha}$ extended to the boundary $\Sigma$, where $\rho$ is any smooth boundary defining function. It is answered from different point of view when a Poincar\'{e}-Einstein manifold is the hyperbolic space \cite{Qi03, ST05, LQS, CLW}. Here we provide a new characterization of this rigidity in terms of the eigenvalues of Dirac operators. If
\begin{equation}\label{HMZ01-1}
\lambda _1 ^2(D) =
\frac{n}{4(n-1)}\mu _1 (L),
\end{equation}
then $(M, \Sigma, g)$ is isometric to the standard hemisphere and $g_+$ is isometric to the hyperbolic space. If
\begin{equation}\label{HMZ02-1}
\kl(\Ds) = \frac{n-1}{2}\frac{Q(M, \Sigma)}{\mbox{Vol}(\Sigma) ^\frac{1}{n-1} },
\end{equation}
then $(M, \Sigma, g)$ is isometric to the flat ball in $\mathbb R^n$ and $g_+$ is isometric to the hyperbolic space.

We point out that the existence of boundary chirality operator on the boundary is weaker than the existence of chirality operator on the whole manifold. Although it is not conformal invariant, the boundary chirality operator yields to a local boundary condition which consists well with the conformal integral Schr\"{o}dinger-Lichnerowicz formula.

The paper is organized as follows. In Section \ref{Preliminaries}, we recall some basic facts about spin manifold, Dirac operator, local boundary condition, integral Schr\"{o}dinger-Lichnerowicz formula and conformal covariance properties of Dirac operator. In Section \ref{PE}, we review the concepts of a conformal compactifiction of a Poincar\'{e}-Einstein manifold and give the proofs of two rigidity results for certain conditions for Ricci curvature and mean curvature. In Section \ref{Main}, we state and prove the main theorems.

\section{Preliminaries}\label{Preliminaries}

In this section, we provide some well-known facts for Dirac operators on manifold with boundary.

\subsection{Dirac operators on manifold with boundary}

Let $(M,g)$ be an $n$-dimensional Riemannian spin manifold with boundary $(\Sigma, \kg)$, where $\kg$ is the induced metric.
Given a spin structure (and so a corresponding orientation) on
manifold $M$, we denote by $\bS _M$ the associated \-spinor bundle,
which is a complex vector bundle of rank $2^{[\frac {n+1}2]}$.
Denote by $\gamma$ the Clifford multiplication
\begin{equation*}
\gamma :\mathbb {C}l(M)\longrightarrow End_{\mathbb C}(\bS _M),
\end{equation*}
which is a fibre preserving
algebra morphism. Let  $\n$ be  the  Riemannian Levi-Civit\`{a} connection of $M$ with respect to the metric $g$ and denote also by the same symbol its corresponding lift to the spinor bundle $\bS_{M}$. It is well known \cite{LM} that there exists a
natural Hermitian metric $\ps{}{}$ on the spinor bundle $\bS_{M}$
which satisfies
\begin{eqnarray}
& X\ps{\psi}{\varphi}= \ps{\n_X\psi}{\varphi}+\ps{\psi}{
	\n_X\varphi}, &\label{comp1}\\
&\langle\gamma(X)\psi,\gamma(X)\varphi\rangle =|X|^2
\ps{\psi}{\varphi},&\label{comp2}\\
&\n_X\big(\gamma(Y)\psi\big)
= \gamma(\n_XY)\psi
+\gamma(Y)\n_X\psi,&\label{comp3}
\end{eqnarray}
for any vector field $X,Y \in \Gamma(TM)$, and for any spinor
fields $\varphi,\psi\in \Gamma(\bS_{M})$. Let $\omega_{n}$ be the complex volume form defined by
\begin{equation}\label{vol-form}
\omega_{n}=\left(\sqrt{-1}\right)^{[\frac{n+1}{2}]}e_1\cdot\ldots \cdot e_n.
\end{equation}
When the dimension $n$ of manifold $M$ is even, the spinor bundle $\bS_{M}$ splits into the direct sum of the subbundles
\begin{equation*}
\bS_{M} =\bS_{M}^+\oplus \bS_{M}^{-},
\end{equation*}
where $\bS_{M}^{\pm}$ are the $\pm1$-eigenspaces of the endomorphism $\gamma(\omega_{n})$.

The Dirac operator $D$ on $\bS_{M}$ is the first order
elliptic differential operator locally given by
\begin{equation*}
D\varphi=\sum_{\kappa =1}^{n}\gamma(e_\kappa){\n}_{\kappa}\varphi
\end{equation*}
for $\varphi\in \Gamma(\bS_{M})$, where $\{e_1,\dots,e_{n}\}$ is a local orthonormal frame of $TM$. When  $n$ is even, the Dirac operator $D$ maps $\bS_{M}^{\pm}$ onto  $\bS_{M}^{\mp}$, i.e. it interchanges positive and negative spinor fields.

The unit normal vector field $e_n$ of hypersurface induces a spin structure on $\Sigma$. Denote the restricted spinor bundle  by $\bS_{\Sigma} =\bS_{M}|_{\Sigma}$. This $\bS _\Sigma$ is referred as the extrinsic spinor bundle of $\Sigma$. We denote also by $\ns$ the spinorial connection acting on the spinor bundle $\bS_{\Sigma}$. The extrinsic spin connection and the extrinsic Dirac operator of $\Sigma$ acting on $\bS _\Sigma$ are given by
\begin{equation}\label{extrinsic-connec}
\nabla ^{\Sigma}:= d +\frac{1}{4} g(\nabla e_i, e_j)\gamma(e_i)\gamma(e_j), 	
\end{equation}
and
\begin{equation}\label{extrinsic-Dirac}
 \Ds= \gamma(e_n) \gamma(e_i) \ns_i.
\end{equation}

As $\Sigma$ equips with the induced spin structure, there is the intrinsic spin bundle $\kS $ on $\Sigma$ with induced spin connection $\kC$ and the Clifford multiplication $\kgamma$. The intrinsic spin connection $\kC$ and the intrinsic Dirac operator $\kD$ of $\Sigma$ acting on $\kS $ are given by
\begin{equation}\label{intrinsic-connec}
\kC:= d +\frac{1}{4} \kg(\kC e_i, e_j) \kgamma(e_i)\kgamma(e_j),
\end{equation}
and
\begin{equation}\label{intrinsic-Dirac}
\kD = \kgamma(e_i)  \kC_i.	
\end{equation}
In general, $\Big(\bS _\Sigma, D^{\Sigma}\Big)$ and $\Big(\kS, \kD \Big)$ are not equivalent. They are isomorphic to each other if $n$ is odd, and the dimension of $\bS _{\Sigma}$ is twice the dimension of $\kS $ if $n$ is even. However, they play the same role. In particular, $\Ds$ and $\kD$ have the same eigenvalues (c.f. \cite{HMZ02}).

The restriction of the spin connection $\nabla$ on $\Sigma$, acting on $\bS _{\Sigma}$, differs with
$\ns$ by the second fundamental forms, i.e., for $\phi \in \Gamma(\bS_{\Sigma})$,
\begin{eqnarray}
\begin{aligned}
\n_i\phi &= \ns_i\phi+\frac{1}{2} g(\nabla _{e_i} e_n, e_j) \gamma(e_n)\gamma(e_j)\phi  \\
         &= \ns_i\phi+\frac12 h_{ij}\gamma(e_n)\gamma(e_j)\phi.\label{Gauss2}
\end{aligned}
\end{eqnarray}
This is called the spinorial Gauss formula. Therefore, on $\Sigma $, for $\phi\in \bS_{\Sigma}$, direct calculation yields
\begin{eqnarray}
\begin{aligned}
 \gamma(e_n)\gamma(e_i)\nabla_i \phi &= \gamma(e_n)\gamma(e_i) \left( \ns_i  + \frac{1}{2}h _{ij} \gamma(e_n)\gamma(e_j) \right)\phi\nonumber\\
                                     &= D ^{\Sigma} \phi -\frac{n-1}{2}H \phi. \label{extrin1}\\
\end{aligned}
\end{eqnarray}
On the other hand,
\begin{eqnarray}
 \begin{aligned}
 \nabla _i (\gamma(e_n) \phi) &=\gamma(\nabla_i e_n)\phi + \gamma(e_n) \nabla_i \phi\nonumber\\
                              &=\gamma(e_n) \n_i \phi +h _{ij} \gamma(e_j) \phi.
 \end{aligned}
\end{eqnarray}
Therefore
\begin{eqnarray}
\begin{aligned}
\ns_i (\gamma(e_n) \phi) &=\n_i (\gamma(e_n) \phi) -\frac{1}{2} g(\n_i e_n, e_j) \gamma(e_n) \gamma(e_j)(\gamma(e_n) \phi)\nonumber\\
&= \gamma(e_n) \n_i \phi -\frac{1}{2}h _{ij} \gamma(e_n) \gamma(e_n) \gamma(e_j) \phi\nonumber\\
 &= \gamma(e_n) \ns_i \phi,
\end{aligned}
\end{eqnarray}
and
\begin{eqnarray}
\begin{aligned}
 D ^ {\Sigma} (\gamma(e_n) \phi)&=\gamma(e_n)\gamma(e_i) \ns_i (\gamma(e_n) \phi)\nonumber\\
 &=\gamma(e_n) \gamma(e_i) \gamma(e_n) \ns_i \phi\nonumber\\
 &=-\gamma(e_n) \Ds \phi.
\end{aligned}
\end{eqnarray}
These yield to the integral Schr\"{o}dinger-Lichnerowicz formula
\begin{eqnarray}
	\int _M  | \nabla \phi |^2d\mu_g =\int _M \left(|D\phi|^2 -\frac{R}{4}|\phi|^2\right)d\mu_g+\int _{\Sigma} \left(\langle \phi, D^ {\Sigma} \phi \rangle- \frac{(n-1)H}{2} |\phi |^2\right)d\sigma_{\gs}.  \label{SL}
\end{eqnarray}

\subsection{Local boundary condition}
It is straightforward to derive
\begin{equation}\label{selfadjoint-non}
\int _M \ps{D\phi}{\psi }d\mu_g - \int _M \ps{\phi}{D\psi}d\mu_g = \int_\Sigma \ps{\gamma(e_n)\phi}{\psi}d\sigma_{\gs}.
\end{equation}
From \eqref{selfadjoint-non}, we know that $D$ is not self-adjoint without posing suitable boundary value. We refer to \cite{APS1, APS2, APS3, GHHP, FS, HMZ01CMP, HMR02, BW, BC, BB} for relevant elliptic boundary conditions.
However, neither the Dirichlet nor the Nermann boundary value makes $D$ elliptic and self-adjoint.

As $D$ is the first order differential operator, and acts on spinors which are vector value functions, the standard
theory of PDEs indicates vanishing of ``half'' vector value
functions on the boundary is elliptic boundary condition. This requires $\bS _\Sigma =\bS ^+ _\Sigma \oplus \bS ^- _\Sigma$,
where $\bS ^\pm _\Sigma$ are two sub spinor bundles of equal dimension.
Then we can take ``half'' part to be zero. This is called local
boundary condition. There is topological obstruction for the existence of local boundary condition to make $D$ self-adjoint.
However, it does exist if the boundary chirality operator exists. An operator $\Gamma \in \text{Hom}(\bS_\Sigma) $
is said to be a boundary chirality operator if it satisfies the
following conditions, for $\phi, \psi \in\bS_\Sigma,$
\begin{equation}\label{def-chi}
\left \{ \aligned
\Gamma ^2 = & Id, \\
\ns_{e_i} \Gamma =&0,
\\
\gamma(e_n) \Gamma  =& - \Gamma \gamma(e_n),\\
\gamma(e_i)
\Gamma  =& \Gamma \gamma(e_i),\\
\ps{\Gamma  \phi}{ \Gamma\psi}
=& \ps{\phi}{\psi}.
\endaligned
\right.
\end{equation}

If the dimension $n$ of $M$ is even, one can always find boundary chirality operator $\Gamma:=\gamma(\omega_{n})\gamma(e_n)$. If $M$ is a spacelike hypersurface with boundary $\Sigma $ and timelike unit normal vector $e_0$ in a Lorentzian manifold. The boundary chirality operator is defined as $\Gamma:=\gamma(e_0) \gamma(e_n)$. In both cases there exists chirality operator globally defined over $M$. However, boundary chirality operator is only defined on the boundary, which is weaker that the existence of chirality operator. Supposing the boundary chirality operator exists, we define
\begin{equation*}
\Gamma_{\pm} ^{loc}=\Big\{\phi \in \bS _{\Sigma}\;:\; \mbox{P} _\pm\phi
=0  \Big\},
\end{equation*}
where $\mbox{P} _\pm$ are pointwise projection operators acting on  $\bS_\Sigma$ defined by
\begin{equation}\label{defP}
\mbox{P} _\pm=\frac{1}{2}\big(\hbox{Id}\mp \Gamma\big).
\end{equation}
It is easy to check that, for $\varphi,\psi\in \Gamma(\bS_\Sigma)$,
\begin{equation}\label{P-orth}
\ps{\mbox{P} _{\pm}\varphi}{\mbox{P} _{\mp}\psi}=0.
\end{equation}
This implies that $\mbox{P} _{+}$ and $\mbox{P} _{-}$ are orthogonal to each other. From \eqref{def-chi}, \eqref{defP} and \eqref{P-orth}, we have
\begin{equation}\label{exchange}
\Ds \mbox{P}_{\pm}=\mbox{P} _{\mp}\,\,\Ds.
\end{equation}
If $\phi \in \Gamma ^{loc}_{\pm}$, then $\gamma(e_n)\phi \in \Gamma ^{loc}
_{\mp}$. Therefore, from \eqref{selfadjoint-non} and \eqref{defP},  $D$ is self-adjoint under the local boundary condition.

It is straightforward that, for $\phi,\psi\in \Gamma(\bS_{\Sigma})$,
 \begin{eqnarray*}
 \begin{aligned}
\langle \gamma(e_i)\gamma(e_j)\phi, \psi\rangle &=-\langle \phi, \gamma(e_i)\gamma(e_j )\psi\rangle, \quad \text{for} \quad i \neq j,\\
	e_i \langle\phi, \psi\rangle &=\langle\n_i \phi, \psi\rangle +\langle\phi, \n_i \psi \rangle \\
                                 &=\langle \ns_i\phi, \psi\rangle +\langle\phi, \ns_i \psi \rangle,\\
\ns_i(\gamma(e_n)\gamma(e_j)\phi)&=\gamma(e_n)\gamma(e_j) \ns_i \phi.
 \end{aligned}
 \end{eqnarray*}
Using \eqref{extrin1}, we can obtain
\begin{equation*}
\begin{aligned}
e_i	\langle  \gamma(e_n)\gamma(e_i) \phi, \psi\rangle = &\langle D ^\Sigma \phi, \psi \rangle -\langle \phi, D ^ \Sigma \psi \rangle\\
		= &\langle \gamma(e_n)\gamma(e_i) \n_i \phi, \psi \rangle -\langle \phi, \gamma(e_n)\gamma(e_i) \n_i \psi \rangle,
	\end{aligned}
\end{equation*}
which imply that  $\Ds $ and $\gamma(e_n)\gamma(e_i)\n_i $ are both self-adjoint on $\Sigma$, i.e., for $\phi,\psi\in \Gamma(\bS_{\Sigma})$,
\begin{equation*}
\int _\Sigma \langle D ^\Sigma \phi, \psi \rangle d\sigma_{\gs} =\int _\Sigma \langle \phi, D^\Sigma \psi \rangle d\sigma_{\gs}, \end{equation*}
and
\begin{equation*}
\int _\Sigma \langle \gamma(e_n)\gamma(e_i)\n_i \phi, \psi \rangle d\sigma_{\gs} =\int _\Sigma \langle \phi, \gamma(e_n)\gamma(e_i)\n_i\psi \rangle d\sigma_{\gs}.	
\end{equation*}

The following theorem is well-known.
\begin{theorem}\label{BVP-spinor}
Suppose $M$ is an $n$-dimensional compact spin manifold with boundary $\Sigma$ which equips with a boundary
chirality operator ($n \geq 3$). Suppose the scalar curvature $R\geq 0$ and the mean
curvature $H \geq 0$. Moreover, either $R>0$ at some point in $M\setminus \Sigma$ or $H>0$ at some point on $\Sigma$. 
Given any $\Phi _0 \in \bS _M$, $\phi _0 \in \bS _{\Sigma}$, there exists a unique smooth spinor $\Psi$ such that\\
\begin{equation*}
\left\{\aligned D\Psi & =\Phi _0 & \mbox{in} \quad M, \\
{\textnormal{P} } _{\pm} \Psi &=\textnormal{P} _{\pm} \phi_0 &\mbox{on}\quad \Sigma.
\endaligned\right.
\end{equation*}
\end{theorem}

\subsection{Conformal covariance of the Dirac operator}
We now recall some properties of  the conformal behavior of spinors on a Riemannian spin manifolds. For more details, we refer to \cite{H74, H86, HMZ01CMP, HMZ02}. Let  $u\in C^\infty(M)$ be a smooth function defined on manifold $M$ and  $\bg=e^{2u}g$ be a conformal change of the metric $g$. This yields the bundle isometry between the two spinor bundles  $\bS_M$ and $\obS_M$, i.e.
\begin{gather*}
\begin{array}{ccc}\label{cc}
\bS_M & \longrightarrow & \obS_M  \\
\varphi & \longmapsto & \overline{\varphi}.
\end{array}
\end{gather*}
We can also relate the corresponding Levi-Civit\`{a}  connections, Clifford multiplications and Hermitian scalar  products.  Denoting by $\oconn$, $\bar{\gamma}$ and $\psb{}{}$ the associated Levi-Civit\`{a}  connection, Clifford multiplication and Hermitian inner product
on sections of the bundle $\obS_M$, one has
\begin{equation*}
\begin{aligned}
\oconn_X\overline{\psi}  =&  \overline{\n_X\psi-\frac{1}{2}\gamma(X) \gamma(\n u)\psi-\frac{1}{2}\ps{X}{\n u}\psi},\\
\bar{\gamma}(\overline{X})\overline{\psi}  =&  \overline{\gamma(X)\psi},\\
\psb{\overline{\psi}}{\overline{\varphi}}  =& \ps{\psi}{\varphi},
\end{aligned}
\end{equation*}
for all $\psi,\varphi\in\Gamma(\bS_M)$, $X\in\Gamma({\rm{T}}M)$ and where $\overline{X}:=e^{-u}X$ denotes the vector field over $(M^n,\bg)$. From these identifications, one has  the relation between the Dirac operators $D$ and $\bar{D}$ acting respectively on sections of~$\bS_M$ and~$\obS_M$, i.e.
\begin{equation}\label{Dirac-rela}
\overline{D}\,\, \overline{\psi}=e^{-\frac{n+1}{2}u}\overline{D\left(e^{\frac{n-1}{2}u}\psi\right)}
\end{equation}
which shows that the Dirac operator is a conformally covariant differential operator.

The conformal change of metric on $M$ induce the corresponding  change of  metric on the hypersurface $\Sigma$, i.e. $\bar{\kg}=e^{2u}\kg$. Denote by $\overline{\Ds}$ the hypersurface Dirac operator acting on the spinor bundle $\obS_\Sigma :=\obS_M|_{\Sigma}$. For the Dirac operators $\Ds$ and $\overline{\Ds}$, we have, for $\psi\in \Gamma(\bS_\Sigma)$,
\begin{equation}\label{conf-Dirac1}
\overline{\Ds}\left(e^{-\frac{n-2}{2}u}\overline{\psi}\right) =e^{-\frac{n}{2}u}\,\,\overline{\Ds\psi},
\end{equation}
which is analogous to \eqref{Dirac-rela}.

Assume that the dimension $n\geq 3$  and $f\in C^\infty(\overline{M})$ is positive function satisfying  $e^{u}=f^\frac{2}{n-2}$. The volume forms of two metrics $\bar{g}$, $g$ and their restriction to the boundary $\Sigma$ satisfy
\begin{equation*}
d\mu_{\bar{g}}= f^{\frac{2n}{n-2}}d\mu_{g}, \quad d\sigma_{\bar{\gs}} = f^{\frac{2(n-1)}{n-2}}d\sigma_{\gs}.
\end{equation*}
The conformal Laplace operator and conformal mean curvature operator obey the conformal transformation laws
\begin{eqnarray}
\bar{L}(f^{-1}v)&=&f^{-\frac{n+2}{n-2}}Lv, \label{conf-Laplace0}\\
\bar{B}(f^{-1}v)&=&f^{-\frac{n}{n-2}}Bv,\label{conf-mOperator0}
\end{eqnarray}
where $v\in C^\infty(M)$.
From \cite{Escobar92}, the scalar curvatures and mean curvature under conformal change  yield
\begin{eqnarray}
\overline{R} &=&f^{-\frac{n+2}{n-2}}Lf,
\label{conf-scalar}\\
\overline{H} &=&f^{-\frac{n}{n-2}}Bf.\label{conf-mean}
\end{eqnarray}
Taking $\psi = f^{-\frac{n-1}{n-2}} \phi$,  by \eqref{Dirac-rela} and \eqref{conf-Dirac1} we have
\begin{equation*}
\overline D\, \overline{\psi} = f^{-\frac{n+1}{n-2}}\overline {D\phi },\quad
\overline{D ^\Sigma } (f^{-1} \overline{\psi}) = f^{-\frac{n}{n-2}} \overline {D^\Sigma \phi }.
\end{equation*}
The Penrose (or twistor) operator $\P$ is defined by
 \begin{equation}\label{Penrose}
 \P_X\phi =\nabla_X\phi  +\frac{1}{n} \gamma(X) D\phi,
 \end{equation}
for any  $X\in \Gamma(TM)$ and $\phi\in \Gamma(\bS_M)$.
The integral Schr\"{o}dinger-Lichnerowicz formula \eqref{SL} can be written as
\begin{equation}\label{SLP}
	\int _{\Sigma} \left(\langle \phi, D^ {\Sigma} \phi \rangle- \frac{(n-1)H}{2} |\phi |^2\right)d\sigma_{\gs} =
\int _M \left(| \P \phi |^2 +\frac{R}{4}|\phi|^2 -\frac{n-1}{n}|D\phi|^2\right)d\mu_{g}.
\end{equation}
Applying \eqref{SLP} to the conformal metric  $\bg$ and $\overline{\psi} \in
\Gamma(\obS_{M})$, it gives
\begin{equation}\label{CSL}
\int _{\Sigma} \left(\langle \overline{\psi}, \overline{D^ {\Sigma}}\, \overline{\psi} \rangle_{\bar{g}}
- \frac{(n-1)}{2} \overline{H}|\overline{\psi} |_{\bar{g}} ^2\right) d\sigma_{\bar{\gs}} =
\int _M \left(| \overline \P \,\overline{\psi} |^2 _{\bg} +\frac{\bar R}{4}|\overline{\psi}|^2 _{\bg} -\frac{n-1}{n}|\overline D\, \overline{\psi}|^2 _{\bg} \right)d\mu_{\bg}
\end{equation}
Since $\psi=f^{-\frac{n-1}{n-2}} \phi =f^{-\frac{1}{n-2}} f^{-1}\phi $, we have
\begin{equation*}
\begin{aligned}
\overline{D^ {\Sigma}} \overline{\phi} &=\bar{\gamma}(\bar{e}_n)\bar{\gamma}(\bar{e}_i) (\bar{e}_i f ^{-\frac{1}{n-2}} )f^{-1} \overline{\phi} + f ^{-\frac{1}{n-2}} \overline{D^ {\Sigma}} (f^{-1}\overline{\phi}) \\
&=-\frac{1}{n-2}f^{-\frac{2n-3}{n-2}} \left(\bar{e}_if\right) \bar{\gamma}(\bar{e}_n)\bar{\gamma}(\bar{e}_{i}) \overline{\phi} + f^{-\frac{1}{n-2}} f^{-\frac{n}{n-2}} \overline {D^ {\Sigma} \phi}.
\end{aligned}
\end{equation*}
Noting that  $\langle \bar \phi, \bar{\gamma}(\bar{e}_n)\bar{\gamma}(\bar{e}_i) \bar \phi \rangle _{\bar{g}}$  is imaginary, we can obtain
\begin{equation*}
\int_\Sigma \langle \bar \psi, \overline{D^ {\Sigma}} \bar \psi \rangle _{\bar{g}}d\sigma_{\bar{\gs}}
=\int _\Sigma f^{-\frac{2n}{n-2}}\langle \bar \phi, \overline{D^ {\Sigma} \phi }\rangle _{\bar{g}} d\sigma_{\bar{\gs}}
=\int _\Sigma f^{-\frac{2}{n-2}}\langle \phi, D^ {\Sigma} \phi  \rangle d\sigma_{\gs}.
\end{equation*}
On the other hand, a direct calculation yields
\begin{equation*}
\begin{aligned}
\overline{H}  |\overline{\psi}| ^2 _{\bar{g}} d\sigma_{\bar{\gs}}&=f^{-\frac{n}{n-2}}f^{-1}B f |\bar \phi|^2 _{\bar{g}} f^{-\frac{2(n-1)}{n-2}}d\sigma_{\bar{\gs}}
=f^{-\frac{2}{n-2}} f^{-1}B f |\phi|^2d\sigma_{\gs},\\
\bar R |\overline{\psi}| ^2 _{\bg}&=f^{-\frac{2(n+1)}{n-2}} f^{-1}L f |\bar \phi|^2 _{\bg} =f^{-\frac{2(n+1)}{n-2}} f^{-1}L f |\phi|^2,\\
|\overline{D}\, \overline{\psi}| ^2 _{\bg} &=f^{-\frac{2(n+1)}{n-2}} |\overline{D \phi}|^2 _{\bg} = f^{-\frac{2(n+1)}{n-2}} |D \phi|^2.
\end{aligned}
\end{equation*}
Finally, we obtain the conformal integral Schr\"{o}dinger-Lichnerowicz formula
\begin{equation}\label{Conf-SL}
\begin{aligned}
\int _{\Sigma} & f^{-\frac{2}{n-2}} \left(\langle \phi, D^ {\Sigma} \phi \rangle- \frac{(n-1)}{2} f^{-1}B f|\phi |^2\right)d\sigma_{\gs} \\
&=
\int _M f^{-\frac{2}{n-2}} \left(| \overline \P \,\overline{\psi}|^2 _{\bg} f ^{\frac{2(n+1)}{n-2}}+\frac{f^{-1}L f}{4}|\phi|^2 -\frac{n-1}{n}|D\phi|^2 \right)d\mu_{g},
\end{aligned}
\end{equation}
where $\psi=f^{-\frac{n-1}{n-2}} \phi$.

\section{Poincare-Einstein metrics and rigidity}\label{PE}

In this section, we study the rigidity for $(M, \Sigma, g)$ as a $C^{3,\alpha}$ conformal compactification of the Poincar\'{e}-Einstein manifolds $(\mathring{M}, g_+)$
under certain curvature assumptions. Denote $\mathring{M}=M \setminus \Sigma$. We assume $(\mathring{M}, g_+)$ is a $n$-dimensional Poincar\'{e}-Einstein manifold ($n\geq 3$):
$$
Ric_{g_+}=-(n-1)g_+\quad \textrm{in $\mathring{M}$},
$$
and $g=\rho^2g_+$ can be $C^{3,\alpha}$ extended to the boundary $\Sigma$ for some smooth boundary defining function $\rho$. Recall $\kg=g|_{\Sigma}$ is denoted as the boundary metric,
 $R^{\Sigma}$ is denoted as the scalar curvature of $(\Sigma, \kg)$ and $E_{ij}$ is denoted as the trace free part of Ricci curvature tensor of $(M,g)$.

\begin{theorem}\label{rigidity1}
If $(M, \Sigma, g)$ is a $C^{3, \alpha}$ conformal compactification of Poincar\'{e}-Einstein manifold $(\mathring{M}, g_+)$ and satisfies
$$
H=0,\quad E:=Ric-\frac{R}{n}g =0,
$$
then $(M, \Sigma, g)$ is isometric to the half sphere $(\mathbb{S}^{n}_+, \mathbb{S}^{n-1}, g_{\mathbb{S}})$ and hence $(\mathring{M}, g_+)$ is isometric to the hyperbolic space $\mathbb{H}^{n}$.
\end{theorem}
\begin{proof}
First by the Gauss-Codazzi equation, $R^{\Sigma}=\frac{n-2}{n}R$ when $H=0$ and $E=0$. Hence $R^{\Sigma}$ is a constant.
Consider the transformation of scalar curvature and  Ricci curvature under conformal change $g=\rho^2g_+$, which gives
\begin{equation}\label{eq.1}
\Delta_{g} \rho =\frac{n}{2} \rho^{-1} \left(|\nabla \rho|^2_{g}-1\right)-\frac{1}{2(n-1)}R\rho,
\end{equation}
\begin{equation}\label{eq.2}
\nabla^2\rho - \frac{1}{n}(\Delta_{g}\rho)g=-\frac{1}{n-2}\rho E=0.
\end{equation}
By identifying a collar neighborhood of $\Sigma$ with $[0,\epsilon)\times \Sigma$, $g$ takes the normal form
$$
g=dr^2+g(r)
$$
where $g(r)$ is a family of metrics on $\Sigma$ with $g(0)=\kg$.
Moreover, according to \cite{Gr1}, $\rho$ has the asymptotical expansion
$$
\rho=r+c_2r^2+c_3r^3+o(r^{3})
$$
where
$$
c_2=-\frac{1}{2(n-1)}H=0, \quad c_3=\frac{1}{6(n-2)}R^{\Sigma}-\frac{1}{6(n-1)}R=-\frac{1}{6n(n-1)}R.
$$
Let
$$
A=\frac{n}{2} \rho^{-1} \left(|\nabla\rho|^2_{g} - 1\right)+ \frac{1}{2(n-1)}R\rho.
$$
Then direct computation shows that
$$
A|_{\Sigma}=0
$$
and
$$
\begin{aligned}
\nabla_i A &=n\rho^{-1}\rho_{ij}\rho^j -\frac{n}{2}  \rho^{-2}\left(|\nabla \rho|^2_{g}-1\right)
\rho_i + \frac{1}{2(n-1)}R\rho_i
\\
&
=\rho^{-1}\left[\Delta_{g}\rho-\frac{n}{2} \rho^{-1}  \left(|\nabla \rho|^2_{g}-1\right)
+\frac{1}{2(n-1)}R\rho
 \right] \rho_i
 =0.
\end{aligned}
$$
Hence $A=A|_{\Sigma}=0$. Thus equations (\ref{eq.1}) and (\ref{eq.2})  become
\begin{equation}\label{eq.3}
\Delta_{g}\rho+\frac{1}{(n-1)}R\rho=0,
\end{equation}
\begin{equation}\label{eq.4}
\nabla^2\rho + \frac{1}{n(n-1)}R\rho g =0.
\end{equation}
Notice that $\rho>0$ in the interior. Hence $R$ must be a positive constant. Up to a constant scaling, we can set $R=n(n-1)$.

Recall that $(M,\Sigma, g)$ is a $C^{3, \alpha}$ compactification of a Poincare-Einstein manifold $(\mathring{M}, g_+)$.
By the boundary regularity theorem given in \cite{CDLS}, $(M, \Sigma,g)$ has umbilic boundary.
Since $H=0$, the boundary is actually totally geodesic.
Take $(\widetilde{M}, \tilde{g})$ to be the double of $(M,g)$ across its boundary and $\tilde{\rho}$ to be the odd extension of $\rho$. Then on $\widetilde{M}$, $\tilde{\rho}$ satisfies the equation
\begin{equation}\label{eq.5}
\tilde{\nabla}^2\rho + \rho\tilde{g} =0.
\end{equation}
This is the standard Obata's equation on closed manifold studied in \cite{Ob1}. Since $\widetilde{M}$ is connected  and $\tilde{\rho}$ is a non-constant solution to (\ref{eq.5}), Obata proved that $(\widetilde{M}, \tilde{g})$  is isometric to the standard sphere
$$
\mathbb{S}^{n}=\{z\in \mathbb{R}^{n+1}: |z|=1\}
$$
 and $\tilde{\rho}$ is the coordinate function $z_1$ up to a rotation and constant scaling. Hence $(M,\Sigma,g)$, which is corresponding to $\tilde{\rho}=z_1\geq 0$, is isometric to the half sphere $(\mathbb{S}^{n}_+, \mathbb{S}^{n-1}, g_{\mathbb{S}})$ and $(\mathring{M}, g_+=\rho^{-2}g)$ is isometric to the standard hyperbolic space $\mathbb{H}^{n}$.

\end{proof}

\begin{theorem}\label{rigidity2}
If $(M, \Sigma, g)$ is a $C^{3, \alpha}$ conformal compactification of Poincar\'{e}-Einstein manifold $(\mathring{M}, g_+)$ and satisfies
$$
H=C,\quad Ric=0,
$$
then $(M, \Sigma, g)$ is isometric to flat ball $(\mathbb{B}^{n}, \mathbb{S}^{n-1}, g_{\mathbb{R}})$ and hence $(\mathring{M}, g_+)$ is isometric to the hyperbolic space $\mathbb{H}^{n}$.
\end{theorem}
\begin{proof}
Notice here $R^{\Sigma}=\frac{n-1}{n}H^2$ by the Gauss-Codazzi equation and hence $R^{\Sigma}$ is a constant.
Consider the transformation of scalar curvature and  Ricci curvature under conformal change $g=\rho^2g_+$, which gives
\begin{equation}\label{eq.7}
\Delta_{g} \rho =\frac{n}{2} \rho^{-1} \left(|\nabla \rho|^2_{g}-1\right),
\end{equation}
\begin{equation}\label{eq.8}
\nabla^2\rho - \frac{1}{n}(\Delta_{g}\rho) g=-\frac{1}{n-2}\rho E=0.
\end{equation}
By identifying a collar neighborhood of $\Sigma$ with $[0,\epsilon)\times \Sigma$, $g$ takes the normal form
\begin{equation}\label{normal}
g=dr^2+g(r)
\end{equation}
where $g(r)$ is a family of metrics on $\Sigma$ with $g(0)=\kg$.
Then according to \cite{Gr1}, $\rho$ has the  asymptotical expansion
$$
\rho=r+c_2r^2+c_3r^3+o(r^{3}),
$$
where
$$
c_2=-\frac{1}{2(n-1)}H, \quad c_3=\frac{1}{6(n-2)}R^{\Sigma}-\frac{1}{6(n-1)}H^2=0.
$$
Direct computation shows that
$$
\Delta_{g}\rho|_{\Sigma}= -\frac{n}{n-1}H,
$$
and
$$
\begin{aligned}
\nabla_i (\Delta_{g}\rho)&=n\rho^{-1}\rho_{ij}\rho^j -\frac{n}{2}  \rho^{-2}\left(|\nabla \rho|^2_{g}-1\right)
\rho_i
\\
&
=\rho^{-1}\left[\Delta_{g}\rho -\frac{n}{2} \rho^{-1}  \left(|\nabla \rho|^2_{g}-1\right)
 \right] \rho_i
 =0.
\end{aligned}
$$
Hence all over $M$,
$$
\Delta_{g}\rho\equiv  -\frac{n}{n-1}H.
$$
Since $\rho>0$ in the interior, we have that $H$ must be a positive constant. Up to a scaling, we can set $H=n-1$ and hence
$\Delta_{g}\rho = -n$. Thus equations (\ref{eq.7}) and (\ref{eq.8}) become
\begin{equation}\label{eq.9}
|\nabla \rho|^2_{g}-1+2\rho=0,
\end{equation}
\begin{equation}\label{eq.10}
\nabla^2\rho +g =0.
\end{equation}
Moreover, $R^{\Sigma}=\frac{n-2}{n-1}H^2=(n-2)(n-1)$ implies that the boundary $(\Sigma,\kg)$ has positive Yamabe constant. By
\cite{WiYa99}, $\Sigma$ is connected.

Take any normal geodesic $\gamma(t)$ such that $\gamma(0)=p\in \Sigma$.
Then $\gamma(t)=(t,p)$. By equation (\ref{eq.10}), the function $f(t)=\rho(\gamma(t))$ satisfies
$$
f''(t)+1=0, \quad f(0)=0,\quad f'(0)=\partial_r\rho|_{\Sigma}= 1.
$$
Hence in the small colloar neighborhood,
\begin{equation}\label{rho}
f(t)=t-\frac{1}{2}t^2, \quad \Rightarrow\quad  \rho=r-\frac{1}{2}r^2.
\end{equation}
On each hypersurfaces $\Sigma_{r}=\{r=constant\}$ for $r$ small, $\rho|_{\Sigma_r}$ is a constant. Moreover, by (\ref{eq.10})  $\rho|_{\Sigma_r}$ satisfies
$$
(\nabla^{\Sigma_r})^2\rho -(\partial_r\rho) h(r)+g(r)=0
$$
where $h(r)$ is the second fundamental form for each level set $(\Sigma_r, g(r))$ w.r.t. outward unit normal $-\partial_r$ and $\nabla^{\Sigma_r}$ is the Levi-Civita connection w.r.t. $(\Sigma_r, g(r))$.
However, we know $h(r)=-\frac{1}{2}g'(r)$ while taking the normal form (\ref{normal}).
This implies that
\begin{equation}\label{ode}
(1-r)g'(r)+2g(r)=0, \quad \Rightarrow\quad g(r)=(1-r)^2 \kg.
\end{equation}
Those formulae (\ref{rho}) and (\ref{ode}) hold in the collar neighborhood such that (\ref{normal}) holds. At any point $0<r_0<1$, if (\ref{ode}) holds, then  (\ref{normal}) extends in a neighborhood $[r_0, r_0+\epsilon)$ and hence  (\ref{rho}) and (\ref{ode}) also can be extended.  The extension will not stop until arriving $r=1$. Therefore,
$$
g=dr^2 + (1-r)^2 \kg, \quad 0\leq r<1.
$$
When $r\rightarrow 1$, $(\Sigma_r,g(r))$ shrink to one point since it is connected,
which corresponds to the unique maximum point of $\rho$. The maximum point is non-degenerate and smooth. Hence
$\kg $ must be the standard sphere metric on $\mathbb{S}^{n-1}$. Therefore, by taking $s=1-r$
$$
(M,g)=([0,1]_s\times \mathbb{S}^{n-1}, g=ds^2+s^2g_{\mathbb{S}})
$$
which is the flat ball of radius one in $\mathbb{R}^{n}$. And $g_+=\rho^{-2}g$ with $\rho=(1-s^2)/2$ shows that
$(\mathring{M}, g_+)$ is the standard hyperbolic space $\mathbb{H}^{n}$.

\end{proof}

\section{Main theorems}\label{Main}

In this section, we firstly re-visit and prove the eigenvalue estimates \eqref{HMZ01} and \eqref{HMZ02} when $\Sigma$ equips with a boundary chirality operator. Then we prove the rigidity of Poincar\'{e}-Einstein manifold when \eqref{HMZ01-1} or \eqref{HMZ02-1} holds.

The following two theorems were proved for $n\geq 3$ with a boundary chirality operator \cite{HMZ01CMP} as well as for $n\geq 2$ with a chirality operator \cite{R06}. Here we provide more accurate statements for $n\geq 2$ and manifolds equip with boundary chirality operators. As boundary chirality operator does not give information of whole manifold as chirality operator does, we can not conclude that manifold is the half sphere when $n\geq 3$ in the equality case \cite{HMR02}.

\begin{theorem}\label{Thm-eigen0}
Let $(M, g)$ be an $n$-dimensional ($n\geq 3$) compact spin manifold with boundary $\Sigma$ which equips with a boundary chirality operator. Suppose that $\mu_1(L)> 0$. Then the first nonzero eigenvalue $\lambda _1(D)$ of the Dirac operator $D$ under the local boundary condition satisfies
\begin{equation}\label{Est-eigen0}
\lambda _1 ^2(D) \geq
\frac{n}{4(n-1)}\mu _1 (L).
\end{equation}
Equality holds if and only if there exists a Killing spinor on $M$ and $\Sigma$ is  minimal.
\end{theorem}

\begin{proof}
The proof follows the main argument in \cite{HMZ01CMP, R06} and we present here for completeness.
For $n\geq 3$, let $f>0$ be the positive solution of \eqref{Eigen-bdry0}. From \eqref{conf-scalar} and \eqref{conf-mean}, we find the scalar and mean curvatures of the conformal metric $\bg=f^{\frac{4}{n-2}}g$ satisfy
\begin{equation*}
\begin{aligned}
\overline{R} =&f^{-\frac{n+2}{n-2}}Lf=\mu_1(L)f^{-\frac{4}{n-2}}>0,\\
\overline{H} =&f^{-\frac{n}{n-2}}Bf=0.
\end{aligned}
\end{equation*}
Now we consider the following eigenvalue problem for Dirac operator with local boundary condition
\begin{equation}\label{Eigenspinor1}
\left\{\aligned D\phi & =\lambda_1(D) \phi & \mbox{in} \quad M, \\
\phi\in  &\Gamma_{\pm} ^{loc}              & \mbox{on} \quad \Sigma.
\endaligned\right.
\end{equation}
Along the boundary $\Sigma$, it is easy to check that $\phi\in \Gamma_{\pm} ^{loc}$ implies $\Ds\phi\in \Gamma_{\mp} ^{loc}$. This gives
\begin{equation*}
 \langle \phi, D^ {\Sigma} \phi \rangle=0.
\end{equation*}
Let $\psi=f^{-\frac{n-1}{n-2}}\phi$. The conformal integral Schr\"{o}dinger-Lichnerowicz formula \eqref{Conf-SL} shows
\begin{equation}\label{Conf-SL1}
\begin{aligned}
0=&\int _{\Sigma}  f^{-\frac{2}{n-2}} \left(\langle \phi, D^ {\Sigma} \phi \rangle- \frac{(n-1)}{2} f^{-1}B f|\phi |^2\right)d\sigma_{\gs} \\
=&
\int _M f^{-\frac{2}{n-2}} \left(| \overline \P \,\overline{\psi}|^2 _{\bg} f ^{\frac{2(n+1)}{n-2}}+\frac{1}{4}f^{-1}L f|\phi|^2 -\frac{n-1}{n}|D\phi|^2 \right)d\mu_{g}\\
\geq &\int _M f^{-\frac{2}{n-2}}\left( \frac{1}{4}f^{-1}L f|\phi|^2 -\frac{n-1}{n}|D\phi|^2 \right)d\mu_{g}\\
=&\int _M f^{-\frac{2}{n-2}}\left( \frac{1}{4}\mu_1(L) -\frac{n-1}{n}\lambda_1^2(D) \right)|\phi|^2d\mu_{g}.
\end{aligned}
\end{equation}
Therefore the inequality holds in \eqref{Est-eigen0}. In the equality case, \eqref{Conf-SL1} gives that
\begin{equation*}
\overline \P_X \,\overline{\psi}=\overline{\nabla}_X\,\overline{\psi}  +\frac{1}{n} \overline\gamma(X) \overline{D}\,\overline{\psi}=0
\end{equation*}
for any $X\in \Gamma(TM)$. Since $\overline{D}\,\overline{\psi}=\lambda_1(D) f^{-\frac{n}{n-2}}\overline{\psi}$, we know that
$\overline{\psi}$ is a Killing spinor. Then the standard argument indicates that $f$ is a constant in $M$ \cite{H86}. Thus
$(M^n,g)$ is Einstein and $\Sigma$ is minimal.

\end{proof}

\begin{theorem}\label{Thm-eigen1}
Let $(M, g)$ be a $2$-dimensional compact oriented surface with boundary $\Sigma$ which equips with a boundary chirality operator. Suppose $\chi(M)>0$. 
Then the first nonzero eigenvalue $\lambda _1(D)$ of the Dirac operator $D$ under the local boundary condition satisfies
\begin{equation}\label{Est-eigen1}
\lambda _1 ^2(D) \geq
\frac{2\pi}{Area(M^2, g)}.
\end{equation}
Equality holds if and only if $(M,\Sigma,g)$ is the half sphere.
\end{theorem}
\begin{proof}
For $n=2$, that conformal changing the metric $\bg=e^{2u}g$ yields the transformation rules for sectional curvature $K$ and geodesic curvature $\kappa$
\begin{equation}\label{conf-2dim}
\left \{ \aligned
e^{2u}\, \overline{K}=&K-\Delta u, \\
e^{u}\, \overline{\kappa} =&\kappa+e_2(u),
\endaligned
\right.
\end{equation}
where $e_2$ is the outer unit normal vector field of $\Sigma$. Let $u$ be the solution of
\begin{equation}\label{solution}
\left \{ \aligned
\Delta u=&K-\frac{1}{Area(M^2, g)}\left(\int_{M}Kd\mu_{g}+\int_{\Sigma}\kappa d\sigma_{\gs}\right),  &\text{in } \quad M \\
e_2(u)=&-\kappa, \quad  &\text{on} \quad   \Sigma,
\endaligned
\right.
\end{equation}
Let $\phi$ be the solution of \eqref{Eigenspinor1} and $\psi=e^{-\frac12 u}\phi$. Applying \eqref{SLP} to the conformal metric $\bg=e^{2u}g$, we obtain
\begin{equation}\label{CSL1}
\int _{\Sigma} \left(e^{u}\langle \overline{\psi}, \overline{D^ {\Sigma}}\, \overline{\psi} \rangle_{\bar{g}}
- \frac{1}{2} e^{u}\overline{\kappa}|\overline{\psi} |_{\bar{g}} ^2\right) d\sigma_{\gs} =
\int _M \left(| \overline \P \,\overline{\psi} |^2 _{\bg} +\frac{\overline{K}}{2}|\overline{\psi}|^2 _{\bg} -\frac{1}{2}|\overline D\, \overline{\psi}|^2 _{\bg} \right)d\mu_{\bg}.
\end{equation}
Since $\langle \phi, D^ {\Sigma} \phi \rangle=0$, $ \ps{\psi}{\gamma(d_\Sigma u)\psi} $ is imaginary and
\begin{equation*}\label{imaginary}
\begin{aligned}
\int _{\Sigma} e^{u}\langle \overline{\psi}, \overline{D^ {\Sigma}}\, \overline{\psi} \rangle_{\bar{g}}d\sigma_{\gs}=&-\frac12\int _{\Sigma}  \ps{\psi}{\gamma(d_\Sigma u)\psi}d\sigma_{\gs}+\int _{\Sigma}e^{-u} \ps{\phi}{D^\Sigma\phi}d\sigma_{\gs},
\end{aligned}
\end{equation*}
we obtain the following identity by taking the real part of \eqref{CSL1}
\begin{equation*}
\begin{aligned}
0=&\int _M \left(| \overline \P \,\overline{\psi} |^2 _{\bg} +\frac{\overline{K}}{2}|\overline{\psi}|^2 _{\bg} -\frac{1}{2}|\overline D\, \overline{\psi}|^2 _{\bg} \right)d\mu_{\bg}\\
\geq & \frac12\int _M\left(\overline{K}e^{2u}-\lambda_1^2(D)\right)\abs{\psi}^2d\mu_{g}.
\end{aligned}
\end{equation*}
By the Gauss-Bonnet formula for surfaces with boundary
$$
\int_{M}Kd\mu_{g}+\int_{\Sigma}\kappa d\sigma_{\gs}=2\pi\chi(M),
$$
we obtain
\begin{equation*}
\begin{aligned}
\frac12\int _M\left(\frac{2\pi\chi(M)}{Area(M^2, g)} -\lambda_1^2(D) \right) \abs{\psi}^2d\mu_{g} \geq 0.
\end{aligned}
\end{equation*}
This gives the second inequality in \eqref{Est-eigen1}.

In the equality case, we deduce that $u$ is constant. Then $K$ is constant and the boundary $\Sigma$ is minimal.
Moreover,  $K=e^{-2u}\overline{K}=e^{-2u}\lambda_1^2(D)>0$.
Consider $(\widetilde{M},\tilde{g})$ being the double of $(M,g)$ across its boundary $\Sigma$. Since $\Sigma$ is minimal and one dimensional, it is
totally geodesic. Thus $(\widetilde{M},\tilde{g})$ is a $C^2$ closed compact manifold which has constant Gaussian curvature $K$.
Therefore, $(\widetilde{M},\tilde{g})$ is isometric to $\mathbb{S}^2$ up to a scaling. Since $\Sigma$ is totally geodesic in $\mathbb{S}^2$, which can only be a great circle. Therefore,  $(M, \Sigma, g)$ is the half sphere.
\end{proof}

The following theorem was proved in \cite{HMZ02} when $\Sigma$ is an internal hypersurface in order to use the unique continuation property of the Dirac operator. Now we prove it when $\Sigma$ is the (usual) boundary of $M$ which the Riemannian structure and spin structure are not necessary products near $\Sigma$.

\begin{theorem}\label{Thm-eigen}
Let $(M,g)$ be an $n$-dimensional ($n \geq 3$) compact spin manifold with boundary $\Sigma$ which equips with a boundary chirality operator. Suppose that $\nu_1(B)> 0$. Then the first positive eigenvalue $\kl(\kD)$ of the intrinsic Dirac operator $\kD$ of $\Sigma $ satisfies
\begin{equation}\label{Est-eigen}
	\kl(\kD) \geq \frac{n-1}{2}\nu _1 (B).
\end{equation}
Equality  implies that $(M,g)$ is conformal to a Ricci flat metric.
\end{theorem}

\begin{proof}
	
The proof follows the main argument in \cite{HMZ02}. Let $f>0$ be the positive solution of \eqref{Eigen-bdry1}.
Let $\bg=f^{\frac{4}{n-2}}g$ be a conformal change of the metric $g$. From \eqref{conf-scalar} and \eqref{conf-mean}, we find its scalar and mean curvatures satisfy
\begin{equation*}
\begin{aligned}
\overline{R} =&f^{-\frac{n+2}{n-2}}Lf=0,\\
\overline{H} =&f^{-\frac{n}{n-2}}Bf=\nu_1(B)f^{-\frac{2}{n-2}}>0.
\end{aligned}
\end{equation*}
Let $\eta = f^{-\frac{1}{n-2}} \phi $, $\psi=f^{-\frac{n-1}{n-2}} \phi$. The conformal integral Schr\"{o}dinger-Lichnerowicz formula \eqref{Conf-SL} reduces to
\begin{equation}\label{Conf-SL2}
\begin{aligned}
\int _{\Sigma} \left(\langle \eta, D^ {\Sigma} \eta \rangle - \frac{n-1}{2}\nu _1 (B) |\eta |^2\right)  d\sigma_{\gs}
= \int _M f^{-\frac{2}{n-2}} \left(| \overline \P \bar \psi |^2 _{\bg} f ^{\frac{2(n+1)}{n-2}}-\frac{n-1}{n}|D\phi|^2 \right)d\mu_{g},
\end{aligned}
\end{equation}
Assume that $\vartheta\in \bS_\Sigma$ is an eigenspinor field associated to $\kl(\Ds)$ over the hypersurface $\Sigma$, i.e. $\Ds\vartheta=\kl(\Ds)\vartheta$. Now we solve the following Dirac equation with local boundary condition
\begin{equation}\label{Equ-Dirac}
\left \{ \begin{aligned}
D \phi = & \,\,0 \qquad \qquad \qquad \, \mbox{in} \quad M, \\
\mbox{P}_+ \phi =& \,\,\mbox{P}_+ (f^{\frac{1}{n-2}} \vartheta) \qquad \mbox{on} \quad \Sigma.
\end{aligned}
\right.	
\end{equation}
The existence of \eqref{Equ-Dirac} follows by showing that $\nu _1 (B)>0$ implies the equation with $\mbox{P}_+ \phi =0$ has trivial solution. Since
$\eta = f^{-\frac{1}{n-2}} \phi $, we have $\mbox{P}_{+}\eta=\mbox{P}_{+}\vartheta$ along the boundary $\Sigma$. From \eqref{exchange}, we have
$\Ds \mbox{P}_\pm\vartheta =\kl(\Ds) \mbox{P}_\mp\vartheta $. From the self-adjointness for $\Ds$, one can get
  \begin{equation}\label{equal}
\kl(\Ds) \int _{\Sigma}|\vartheta ^+|^2 =\kl(\Ds) \int _{\Sigma}|\vartheta ^-|^2.
  \end{equation}
By the Cauchy-Schwartz inequality, we have
\begin{equation}\label{CS1}
\begin{aligned}
\int_{\Sigma}\ps{\Ds\eta}{\eta}d\sigma=&2\Re\int_{\Sigma}\ps{\Ds \mbox{P}_{+}\vartheta}{\mbox{P}_{-}\eta}d\sigma_{\gs}\\
=&2\kl\Re\int_{\Sigma}\ps{ \mbox{P}_{-}\vartheta}{\mbox{P}_{-}\eta}d\sigma_{\gs}\\
\leq & \kl(\Ds) \int_{\Sigma}\left(\abs{\mbox{P}_{-}\vartheta}^2+\abs{\mbox{P}_{-}\eta}^2\right)d\sigma_{\gs}\\
=&\kl(\Ds) \int_{\Sigma}\left(\abs{\mbox{P}_{+}\vartheta}^2+\abs{\mbox{P}_{-}\eta}^2\right)d\sigma_{\gs}\\
=&\kl(\Ds) \int_{\Sigma}\abs{\eta}^2d\sigma_{\gs},
\end{aligned}
\end{equation}
Now \eqref{Eigen-bdry1}, \eqref{Conf-SL2}, \eqref{Equ-Dirac} and  \eqref{CS1} indicate that
\begin{equation}\label{inequ0}
\begin{aligned}
0\leq &\int _M | \overline \P \,\overline{\psi} |^2 _{\bg} d\mu_{\bg}\\
\leq &\int _{\Sigma} \left(\langle D^ {\Sigma} \eta,  \eta \rangle - \frac{n-1}{2}\nu _1 (B) |\eta |^2\right)  d\sigma_{\gs}\\
\leq & \int_{\Sigma}\left(\kl(\Ds)- \frac{n-1}{2}\nu _1 (B) \right)|\eta |^2 d\sigma_{\gs}.
\end{aligned}
\end{equation}
Since $\kl(\kD)=\kl(\Ds)$, \eqref{Est-eigen} follows. In the equality case, $\overline{\psi}$ is a parallel spinor field with respect to the  conformal metric $\bg$. Hence $(M,\bg)$ is Ricci flat.
\end{proof}

Now we prove the following two rigidity theorems for Poincar\'{e}-Einstein manifolds.

\begin{theorem}\label{main-rigidity1}
Let $(M,g)$ be an $n$-dimensional ($n \geq 3$) compact spin manifold with boundary $\Sigma$ which equips with a boundary chirality operator.
If $(M, \Sigma, g)$ is a $C^{3, \alpha}$ conformal compactification of Poincar\'{e}-Einstein manifold $(\mathring{M}, g_+)$ and satisfies
\begin{equation*}
\lambda _1 ^2(D) =
\frac{n}{4(n-1)}\mu _1 (L),
\end{equation*}
then $(M, \Sigma, g)$ is isometric to the half sphere $(\mathbb{S}^{n}_+, \mathbb{S}^{n-1}, g_{\mathbb{S}})$ and hence $(\mathring{M}, g_+)$ is isometric to the hyperbolic space $\mathbb{H}^{n}$.
\end{theorem}
\begin{proof}
It is known from Theorem \ref{Thm-eigen0} that $M$ is Einstein and $\Sigma$ is minimal. Then the theorem follows from Theorem \ref{rigidity1}.
\end{proof}

\begin{theorem}\label{main-rigidity2}
Let $(M,g)$ be an $n$-dimensional ($n \geq 3$) compact spin manifold with boundary $\Sigma$ which equips with a boundary chirality operator.
If $(M, \Sigma, g)$ is a $C^{3, \alpha}$ conformal compactification of Poincar\'{e}-Einstein manifold $(\mathring{M}, g_+)$ and satisfies
\begin{equation*}
\kl(\Ds) = \frac{n-1}{2}\frac{Q(M, \Sigma)}{\mbox{Vol}(\Sigma) ^\frac{1}{n-1} },
\end{equation*}
then $(M, \Sigma, g)$ is isometric to flat ball $(\mathbb{B}^{n}, \mathbb{S}^{n-1}, g_{\mathbb{R}})$ and hence $(\mathring{M}, g_+)$ is isometric to the hyperbolic space $\mathbb{H}^{n}$.
\end{theorem}
\begin{proof}
The equality implies that
\begin{equation*}
\kl(\Ds) = \frac{n-1}{2}\nu _1 (B)=\frac{n-1}{2}\frac{Q(M, \Sigma)}{\mbox{Vol}(\Sigma) ^\frac{1}{n-1} }.
\end{equation*}
Thus, from the first equality and Theorem \ref{Thm-eigen}, we know that $\bar g$ is Ricci flat. The second equality implies that $f$ is constant on $\Sigma$, hence $\bar H$ is constant. Therefore the theorem follows from Theorem \ref{rigidity2}.
\end{proof}

{\small {\it Acknowledgement.}
The work of D. Chen was supported by NSF of China grant 11471180 and 11831005.	
The work of F. Wang was supported by NSF of China grant 11571233 and 11871331.
The work of X. Zhang was supported by NSF of China grants 11571345, 11731001 and HLM, NCMIS, CEMS, HCMS of Chinese Academy of Sciences.}

\end{document}